\documentclass[reqno]{amsart}
\usepackage{amsmath,amssymb,amsthm,amscd,amsfonts,mathrsfs,verbatim}
\usepackage[all]{xy}

\newcommand{\Z}{{\mathbb Z}}                   
\newcommand{\R}{{\mathbb R}}                   
\newcommand{\C}{{\mathbb C}}                   
\newcommand{\Mod}{{\mathcal M}}               
\renewcommand{\O}{{\mathcal O}}
\newcommand{\CP}[1]{\mathbb{C}P^{#1}}      
\newcommand{\E}{{\mathcal E}}

\DeclareMathOperator{\Dol}{Dol}
\DeclareMathOperator{\Betti}{B}

\newtheorem{lem}{Lemma}
\newtheorem{theorem}{Theorem}

\title[Geometric P=W conjecture via abelianization]{Simpson's geometric P=W conjecture in the Painlev\'e VI case via abelianization}

\author{Szil\'ard Szab\'o}
\address{Budapest University of Technology and Economics, 1111. Budapest,
Egry J\'ozsef utca 1. H \'ep\"ulet, Hungary, and 
R\'enyi Institute of Mathematics, 1053. Budapest, Re\'altanoda
utca 13-15. Hungary}
\email{szabosz@math.bme.hu, szabo.szilard@renyi.mta.hu}

\begin{document}

\begin{abstract}
We use abelianization of Higgs bundles near infinity to prove the homotopy commutativity assertion of Simpson's geometric P=W conjecture in the Painlev\'e VI case. 
\end{abstract}

\maketitle

\section{Introduction and statement of the main result}

Throughout this paper we will be dealing with the Dolbeault moduli space $\Mod_{\Dol}$ and the character variety (or Betti moduli space) ${\Mod}_{\Betti}$ associated to the 
Painlev\'e VI equation. We fix some point $t \in \CP1 \setminus \{ 0,1,\infty \}$. 
The  Dolbeault moduli space $\Mod_{\Dol}$ parameterizes certain \emph{logarithmic Higgs bundles} $(\E, \theta )$ of rank $2$ over $\CP1$, with logarithmic divisor $0 + 1 + t + \infty$ 
up to isomorphism (or more generally, S-equivalence). 
${\Mod}_{\Betti}$ parameterizes certain representations 
$$
\rho\colon \pi_1 (\CP1 \setminus \{ 0 , 1 , t , \infty \}) \to \mbox{GL}(2,\C ), 
$$
up to overall conjugation. 
The definition of these spaces depends on choices of eigenvalues of the residues of the Higgs field (equivalently, eigenvalues of the 
local monodromy matrices) at the punctures; we choose these eigenvalues of the residues to be all equal to $0$, so that the residues are nilpotent. 
The Dolbeault space also depends on choices of parabolic weights; we choose these weights so that the moduli space be a smooth manifold. 
To these data there correspond eigenvalues of the local monodromy matrices and weights on the Betti side by some transformation~\cite{Sim_Hodge} that we omit to spell out. 
It is known that then these spaces are $\C$-analytic surfaces diffeomorphic to each other under non-abelian Hodge theory composed by the Riemann--Hilbert correspondence: 
$$
  \psi: \Mod_{\Dol} \to {\Mod}_{\Betti}; 
$$
notice however, that $\psi$ is not $\C$-analytic. 
Furthermore, it is known that there exists a proper map 
\begin{align}
   h: \Mod_{\Dol} & \to \C \label{eq:Hitchin} \\ 
   (\E, \theta ) & \mapsto \det (\theta ) \notag
\end{align}
called \emph{Hitchin map}, with elliptic curves as generic fibers. For details, see Section~\ref{sec:proof}. 

For $R \in \R$ let us denote by $B_R(0)$ the open disk of radius $R$ centered at $0$. 
For a simplicial complex $\mathcal{N}$, we denote by $| \mathcal{N} |$ its body, i.e. the associated topological space. 
It essentially goes back to Fricke and Klein~\cite{FK} that there exists a smooth compactification $\widetilde{\Mod}_{\Betti}$ of 
$\Mod_{\Betti}$ by a simple normal crossing divisor $\tilde{D}_{\infty}$ with nerve complex $\mathcal{N}$ such that 
$| \mathcal{N} |$ is homotopy equivalent to $S^1$. 
The $0$-skeleton $\mathcal{N}_0$ is in bijection with irreducible components of $\tilde{D}_{\infty}$, the edges in the $1$-skeleton $\mathcal{N}_1$ 
corresponds to intersection points between the corresponding components. 
Let $T_j \subset {\Mod}_{\Betti}$ denote a sufficiently small punctured tubular neighbourhood of the $j$'th irreducible component of $\tilde{D}_{\infty}$ and set 
$$
  T = \cup_{j \in \mathcal{N}_0} T_j, 
$$
a punctured tubular neighbourhood of $\tilde{D}_{\infty}$. 
Let $\{ \phi_j \}_{j \in \mathcal{N}_0}$ be a partition of unity on $T$ dominated by the covering $\{ T_j \}_{j \in \mathcal{N}_0}$. 
In~\cite{Sim}, in the more general setup of character varieties with an arbitrary number of punctures on the projective line, C.\;Simpson defines the continuous map 
$$
  \phi = (\phi_j)_{j \in \mathcal{N}_0}: T \to \R^{\mathcal{N}_0} 
$$
whose image is the body of a simplicial complex isomorphic to $\mathcal{N}$. 

\begin{theorem}\label{thm:Simpson}
 For all sufficiently large $R \in \R$, there exists a homotopy commutative square  
 $$
  \xymatrix{\Mod_{\Dol} \setminus h^{-1} (B_R(0))  \ar[d]_h \ar[r]^{\psi} & {\Mod}_{\Betti} \setminus \psi( h^{-1} (B_R(0))) \ar[d]^{\phi} \\
  \C \setminus B_R(0) \ar[r] & | \mathcal{N} |. }
 $$
\end{theorem}

This statement in higher generality was conjectured by L.\;Katzarkov, A.\;Noll, P.\;Pandit and C.\;Simpson~\cite[Conjecture 1.1]{KNPS} 
and subsequently named Geometric $P = W$ conjecture by C.\;Simpson~\cite[Conjecture 11.1]{Sim}. 
In~\cite{Sz_PW} the author proved this homotopy commutativity assertion for all Painlev\'e cases, which includes in particular the Painlev\'e VI case treated here. 
The proof in~\cite{Sz_PW} relies on an explicit understanding of the spaces on both sides in terms of low-dimensional topology. 
Here, we offer a novel approach to this result using abelianization, that may lend itself more easily to generalization. 

\noindent {\bf Acknowledgements:} The author would like to thank R.\;Mazzeo, T.\;Mochizuki and C.\;Simpson for useful discussions. 
 During the preparation of this manuscript, the author benefited of support by the \emph{Lend\"ulet} Low Dimensional Topology 
  grant of the Hungarian Academy of Sciences and by the grants K120697 and KKP126683 of NKFIH.

\section{Proof of Theorem~\ref{thm:Simpson}}\label{sec:proof}

Consider $X = \CP1$ with charts $z$ and $w = z^{-1}$, and let $\E$ denote a holomorphic vector bundle of degree $0$ over $\CP1$. 
We denote by $\O$ and $K$ the sheaves of holomorphic functions and holomorphic $1$-forms respectively on $\CP1$. 
We consider parabolically stable logarithmic Higgs fields $\theta$ on $\E$ with singularities at $0,1,t,\infty$ having nilpotent residue at all these points. 
We have 
$$
  \mbox{tr} (\theta ) \in H^0 (\CP1, K  ) = 0, 
$$
and 
$$
  \det (\theta ) \in H^0 (\CP1, K^2  (0 + 1 + t + \infty ) ) \cong \C.
$$
We call this affine space the \emph{Hitchin base}. We fix an isomorphism 
$$
  \O  \cong  K^2  (0 + 1 + t + \infty )
$$
given by 
$$
  1 \mapsto \frac{(\mbox{d} z)^{\otimes 2}}{z(z-1)(z-t)}. 
$$
On the other hand, we consider the holomorphic line bundle $L = K (0 + 1 + t + \infty )$ with the natural projection 
$$
  p_L : \mbox{Tot} ( L ) \to \CP1 
$$
of its total space $\mbox{Tot} ( L )$ to $\CP1$ and denote by 
$$
  \zeta \frac{\mbox{d} z}{z(z-1)(z-t)}
$$ 
the canonical section of $p_L^* L$ over $p_L^{-1} (\C)$. Then the curve 
$$
   \tilde{X}  = \{ (z, \zeta): \quad \zeta^2 + z(z-1)(z-t) = 0 \} \subset \C_z \times \C_{\zeta} 
$$
has a smooth compactification in $\mbox{Tot} ( L )$ at $z = \infty$ by the point 
$$
  w = 0, \zeta = 0. 
$$
We continue to denote this compactification by $\tilde{X}$ and moreover denote by 
\begin{align}\label{eq:double_cover}
   p: \tilde{X} & \to \CP1 \\
  (z, \zeta) & \mapsto z \notag
\end{align}
the restriction of $p_L$, a ramified double covering map with branch points $\{ 0, 1, t, \infty \}$. 
Topologically, $\tilde{X}$ is diffeomorphic to a $2$-torus. 
Let us introduce the bivalued holomorphic $1$-form on $\C \setminus \{ 0,1,t \}$ 
$$
  \omega = \frac{\mbox{d} z}{\sqrt{z(z-1)(z-t)}}. 
$$
Notice that 
$$
  p^* \omega = \frac{\mbox{d} z}{\zeta}
$$ 
is a univalued meromorphic $1$-form on $\tilde{X}$, with simple poles at the branch points of $p$. 

For $R>>0, \varphi \in \R / 2 \pi \Z$ we let $(\E_{R, \varphi} ,\theta_{R, \varphi})$ be any rank $2$ logarithmic Higgs bundle over $\CP1$ with 
$$
  \det (\theta_{R, \varphi} ) = - R e^{\sqrt{-1} \varphi} \in H^0 (\CP1, K^2  (0 + 1 + t + \infty ) )
$$
(the sign is introduced for convenience). 
We will fix $R$ and let $\varphi$ vary, and we assume that $\theta_{R, \varphi}$ depends smoothly and 
$2\pi$-periodically on $\varphi$, thus providing a smooth section of~\eqref{eq:Hitchin} over $|z| = R$; 
such lifts clearly exist. 
The \emph{spectral curve} of $(\E_{R, \varphi} ,\theta_{R, \varphi})$ defined as 
$$
  \tilde{X}_{R, \varphi}  = \left\{ (z, \zeta): \quad \det \left( \theta_{R, \varphi} - \zeta \frac{\mbox{d} z}{z(z-1)(z-t)} \right) = 0 \right\} \subset \mbox{Tot} ( L )
$$
is obtained by rescaling~\eqref{eq:double_cover} in the $\zeta$-direction by the factor $\sqrt{R} e^{\sqrt{-1} \varphi / 2}$. 
In particular, for any $(R, \varphi)$ the branch points of $\tilde{X}_{R, \varphi}$ are the distinct points $\{ 0, 1, t, \infty \}$, and 
the curve $\tilde{X}_{R, \varphi}$ is smooth. 
Let $z_0 \notin \{ 0,1,t, \infty \}$ and fix $\varepsilon_0 > 0$ such that $B_{2\varepsilon_0 } (z_0)$ is disjoint from $\{ 0,1,t, \infty \}$. 
According to \cite[Theorem~1.4]{Moc}, for $z\in B_{\varepsilon_0 } (z_0)$ there exists a smoothly varying frame $e_1 (z), e_2 (z)$ of $\E$ with respect to which we have 
the asymptotic equality 
$$
  \theta_{R, \varphi} (z) - \begin{pmatrix}
                             \sqrt{R} e^{\sqrt{-1} \varphi / 2} & 0 \\
                             0 & - \sqrt{R} e^{\sqrt{-1} \varphi / 2}
                            \end{pmatrix}
	\omega \to 0 
$$
as $R\to\infty$, with exponential rate. 
Because $\omega$ is bivalued, the vectors $e_1 (z), e_2 (z)$ get interchanged as the position of the point $z_0$ moves along a simple loop $\gamma$ around one of the punctures. 
In different terms the monodromy transformation of the trivialization along any such $\gamma$ is the transposition matrix 
$$
  T = \begin{pmatrix}
   0 & 1 \\
   1 & 0
  \end{pmatrix}. 
$$

We consider the connection associated by non-abelian Hodge theory to $(\E_{R, \varphi} ,\theta_{R, \varphi})$ with respect to the gauge $e_1 (z), e_2 (z)$. 
Its connection form is 
\begin{align*}
   a_{R, \varphi} (z, \bar{z}) & = \theta_{R, \varphi} (z) + \overline{\theta_{R, \varphi} (z)} + b_{R, \varphi} \\
   & \approx \sqrt{R}  
   \begin{pmatrix}
    e^{\sqrt{-1} \varphi / 2} \omega + e^{-\sqrt{-1} \varphi / 2} \bar{\omega} & 0 \\
    0 & - e^{\sqrt{-1} \varphi / 2} \omega - e^{-\sqrt{-1} \varphi / 2} \bar{\omega}
   \end{pmatrix} 
   + b_{R, \varphi} 
\end{align*}
where $\approx$ means that the difference of the two sides converges exponentially to $0$ as $R\to\infty$, 
and $b_{R, \varphi}$ is the connection form of the Chern-connection $\partial_{R, \varphi}$ of $(\E_{R, \varphi}, h_{R, \varphi} )$. 
It also follows from \cite[Theorem~1.4]{Moc} that 
\begin{equation}\label{eq:asymptotic_commuting}
   \left[ \begin{pmatrix}
    e^{\sqrt{-1} \varphi / 2} \omega + e^{-\sqrt{-1} \varphi / 2} \bar{\omega} & 0 \\
    0 & - e^{\sqrt{-1} \varphi / 2} \omega - e^{-\sqrt{-1} \varphi / 2} \bar{\omega}
   \end{pmatrix} ,
   b_{R, \varphi} \right] \to 0 
\end{equation}
exponentially as $R\to\infty$. 
Furthermore, as $\mbox{tr} (\theta ) \equiv 0$ the Higgs field induced by $\theta_{R, \varphi}$ on $\det (\E_{R, \varphi} )$ is identically zero. 
It follows that the corresponding Hermitian--Einstein metric $h_{\det (\E )} \equiv 1$  and thus $b_{R, \varphi}$ takes values in $\mathfrak{su}(2, \C)$. 
As we will only be interested in the absolute value of the integral of $a_{R, \varphi}$ along loops and the monodromy of the Chern connection is unitary, 
we will see that the actual shape of $b_{R, \varphi}$ is irrelevant for our purposes. 

In order to get a hold on $\psi (\E_{R, \varphi} ,\theta_{R, \varphi})$, we need to apply the Riemann--Hilbert correspondence to the connection obtained in the previous paragraph. 
For this purpose, we now fix $z_0$ and simple loops $\gamma_0, \gamma_1, \gamma_t$ based at $z_0$ winding about the punctures $\{ 0,1,t \}$ respectively, in positive direction. 
The monodromy matrices of the connection $\mbox{d} + a_{R, \varphi}$ associated to the punctures are given by 
\begin{equation*}
  B_j(R, \varphi) = \exp \oint_{\gamma_j} - a_{R, \varphi} (z, \bar{z})
\end{equation*}
where $j\in \{ 0,1,t \}$. Let us introduce the \emph{half-period integrals} 
$$
  \pi_j  = \oint_{\gamma_j} \omega. 
$$
By this we mean that we fix any one of the two lifts $\tilde{\gamma}_j$ of $\gamma_j$ by~\eqref{eq:double_cover} and set 
$$
  \pi_j = \int_{\tilde{\gamma}_j} p^* \omega. 
$$
By Baker--Campbell--Hausdorff formula and~\eqref{eq:asymptotic_commuting}, as $R\to\infty$ the monodromy matrix $B_j (R, \varphi)$ is asymptotically equal to 
\begin{equation*}
  T A_j(R, \varphi) \exp \sqrt{R} 
  \begin{pmatrix}
   -  e^{\sqrt{-1} \varphi / 2} \pi_j - e^{-\sqrt{-1} \varphi / 2}  \overline{\pi_j} & 0 \\
   0 & e^{\sqrt{-1} \varphi / 2} \pi_j + e^{-\sqrt{-1} \varphi / 2}  \overline{\pi_j}
  \end{pmatrix} 
\end{equation*}
for some special unitary matrix 
\begin{equation}\label{eq:unitary_factor}
  A_j(R, \varphi) = \exp \oint_{\gamma_j} - b_{R, \varphi} (z, \bar{z}) \in \mbox{SU} (2, \C)
\end{equation}
commuting with the matrix on its right. These properties show that necessarily 
$$
   A_j(R, \varphi) = \begin{pmatrix}
                      e^{\sqrt{-1}\mu_j} & 0 \\
                      0 & e^{-\sqrt{-1}\mu_j}
                     \end{pmatrix}
$$
for some $\mu_j = \mu_j (R, \varphi) \in \R$. 
We then get 
\begin{equation*}
 B_j (R, \varphi) \approx 
 \begin{pmatrix}
       0 & \exp \left( - \sqrt{-1}\mu_j + 2\sqrt{R}  \Re ( e^{\sqrt{-1} \varphi / 2} \pi_j) \right) \\
       \exp \left( \sqrt{-1}\mu_j - 2\sqrt{R}  \Re ( e^{\sqrt{-1} \varphi / 2} \pi_j) \right) & 0
      \end{pmatrix} , 
\end{equation*}
and it follows that 
$$
  B_0 (R, \varphi) B_1 (R, \varphi) \approx 
  \begin{pmatrix}
   d_{01} (R, \varphi) & 0 \\
   0 & \frac 1{d_{01} (R, \varphi)} 
  \end{pmatrix}. 
$$
with 
$$
  d_{01} (R, \varphi) = \exp \left( \sqrt{-1} (\mu_1 - \mu_0) + 2\sqrt{R}  \Re ( e^{\sqrt{-1} \varphi / 2} (\pi_0 - \pi_1)) \right).
$$
Therefore, setting 
$$
  X_1 (R, \varphi) = \mbox{tr} (B_0 (R, \varphi) B_1 (R, \varphi) ) 
$$
we find 
\begin{align}
  X_1 (R, \varphi) & \approx 2 \cosh d_{01} (R, \varphi) \notag \\
   & = 2 \cosh \left( \sqrt{-1} (\mu_1 - \mu_0) + 2\sqrt{R}  \Re ( e^{\sqrt{-1} \varphi / 2} (\pi_0 - \pi_1)) \right) .\label{eq:X1}
\end{align}
Similarly, we find 
\begin{align}
 X_2 (R, \varphi) & = \mbox{tr} (B_t (R, \varphi) B_0 (R, \varphi) ) \notag \\ 
  & \approx 2 \cosh \left( \sqrt{-1} (\mu_0 - \mu_t) + 2\sqrt{R} \Re (e^{\sqrt{-1} \varphi / 2} (\pi_t - \pi_0 )) \right) \label{eq:X2} \\
 X_3 (R, \varphi) & = \mbox{tr} (B_1 (R, \varphi) B_t (R, \varphi) ) \notag \\ 
 & \approx 2 \cosh \left( \sqrt{-1} (\mu_t - \mu_1) + 2\sqrt{R} \Re ( e^{\sqrt{-1} \varphi / 2}(\pi_1 - \pi_t)) \right) \label{eq:X3}
\end{align}
where $X_2 (R, \varphi)$ and $X_3 (R, \varphi)$ are defined by the equalities in these formulas. 
It is known from \cite{FK} that these quantities fulfill the equation 
$$
  X_1 X_2 X_3 + X_1^2 + X_2^2 + X_3^2 - s_1 X_1 - s_2 X_2 - s_3 X_3 + s_4 = 0
$$
for some constants $s_1, s_2, s_3, s_4 \in \C$. Let us denote by 
$$
  \mathcal{M} = \{ (X_1, X_2, X_3) : \quad   X_1 X_2 X_3 + X_1^2 + X_2^2 + X_3^2 - s_1 X_1 - s_2 X_2 - s_3 X_3 + s_4 = 0 \} \subset \C^3 
$$ 
the affine surface determined by this equation. The natural projectivization $\overline{\mathcal{M}}$ of $\mathcal{M}$ in $\CP3$ is given by the equation
$$
    X_1 X_2 X_3 + X_0 X_1^2 + X_0 X_2^2 + X_0 X_3^2 - s_1 X_0^2 X_1 - s_2 X_0^2 X_2 - s_3 X_0^2 X_3 + s_4 X_0^3 = 0 
$$
where $[X_0 : X_1 : X_2 : X_3] \in \CP3$. 
The projective surface $\overline{\mathcal{M}}$ contains $\mathcal{M}$ as the points with $X_0 \neq 0$, and we may then suppose $X_0 = 1$. 
A sequence of points of $\mathcal{M}$ with $X_1\to\infty$ (or $X_2\to\infty$ or $X_3\to\infty$) may converge to a point of 
$$
  \bar{D}_{\infty} = \overline{\mathcal{M}} \setminus \mathcal{M}. 
$$
We call $\bar{D}_{\infty}$ the \emph{divisor at infinity} of $\overline{\mathcal{M}}$. Clearly, 
$$
  \bar{D}_{\infty} = \{ [0: X_1 : X_2 : X_3]:  \quad X_1 X_2 X_3 = 0 \} \subset \CP2_{\infty} \subset \CP3
$$
under the map 
$$
  \CP2_{\infty} \ni [X_1 : X_2 : X_3] \mapsto [0: X_1 : X_2 : X_3] \in \CP3. 
$$
In order to obtain a smooth compactification 
$$
  \mathcal{M} \subset \widetilde{\mathcal{M}}, 
$$
in principle one may need to apply blow-up of $\overline{\mathcal{M}}$ at points of $\mathcal{M}$ and some of the points 
$$
  [0:1:0:0], \quad [0:0:1:0], \quad [0:0:0:1]. 
$$
As it can be checked \cite{Sz_PW}, in the Painlev\'e VI case no such blow-ups are needed and we have $\widetilde{\mathcal{M}} = \overline{\mathcal{M}}$. 
Let $\mathcal{N} (\bar{D}_{\infty})$ stand for the dual simplicial complex of $\bar{D}_{\infty}$. The vertices of $\mathcal{N} (\bar{D}_{\infty})$ are given by 
$$
   v_1 = [0 : 0: X_2 : X_3], \quad v_2 = [0 : X_1 : 0 : X_3], \quad v_3 = [0:X_1 : X_2 :0],
$$
where for instance $[0:X_1 : X_2 :0]$ means the corresponding line in $\CP2_{\infty}$. Let the edges of $\mathcal{N} (\bar{D}_{\infty})$ be denoted by 
$$
  [v_1 v_2], \quad [v_2 v_3], \quad [v_3 v_1], 
$$
respectively corresponding to the following intersection points of the divisor components: 
$$
  [0:0:0:1], \quad [0:1:0:0], \quad   [0:0:1:0]  . 
$$
For $1 \leq j \leq 3$, let $\phi_j$ be continuous maps from a neighbourhood $T_j$ of the line encoded by $v_j$ into $[0,1]$, satisfying 
$$
  \phi_1 + \phi_2 + \phi_3 = 1 
$$
and set 
$$
  \phi = \begin{pmatrix}
          \phi_1 \\ \phi_2 \\ \phi_3
         \end{pmatrix}.
$$

We need to show that for any section $(\E_{R, \varphi} ,\theta_{R, \varphi})$ of $h$ over $|z| = R$, the loop 
$$
  \phi \circ \psi (\E_{R, \varphi} ,\theta_{R, \varphi})
$$
is a generator of $\pi_1 (| \mathcal{N} (\bar{D}_{\infty}) |)$. 
For this purpose, we need to study the asymptotic behaviour of the quotients 
$$
  \frac{X_1(R, \varphi)}{X_2(R, \varphi)}, \quad \frac{X_3(R, \varphi)}{X_1(R, \varphi)}, \quad \frac{X_2(R, \varphi)}{X_3(R, \varphi)}
$$
as $R\to +\infty$, and in particular the way this behaviour depends on $\varphi \in [0, 2\pi ]$. 
Notice that for $d\in \C$ with  
$|\Re (d)| >> 0$ we have 
$$
  |2 \cosh (d ) | \approx e^{| d |}.
$$
Applying this asymptotic equivalence to~\eqref{eq:X1},~\eqref{eq:X2},~\eqref{eq:X3} for $R>>0$ we find 
\begin{align*}
 |X_1(R, \varphi)| & \approx  \exp \left( 2\sqrt{R} |\Re ( e^{\sqrt{-1} \varphi / 2} ( \pi_0 - \pi_1))|\right), \\
 |X_2(R, \varphi)| & \approx  \exp \left( 2\sqrt{R} |\Re ( e^{\sqrt{-1} \varphi / 2} ( \pi_t - \pi_0))|\right), \\
 |X_3(R, \varphi)| & \approx  \exp \left( 2\sqrt{R} |\Re ( e^{\sqrt{-1} \varphi / 2} ( \pi_1 - \pi_t))|\right).  
\end{align*}
For generic $t\in \CP1 \setminus \{ 0, 1 , \infty \}$ the periods $\pi_0, \pi_1, \pi_t$ are not colinear in $\C$, 
said differently they form the vertices of a non-degenerate triangle $\Delta$ with sides 
$$
  a = \pi_0 - \pi_1, \quad   b = \pi_t - \pi_0, \quad c = \pi_1 - \pi_t.
$$
Consider the triangles $e^{\sqrt{-1} \varphi / 2} \Delta$ as $\varphi$ ranges over $[0, 2\pi )$. 
A straightforward geometric inspection shows that the lengths of the projection onto the real axis 
of the three sides of $e^{\sqrt{-1} \varphi / 2} \Delta$ obey the following rule. 
\begin{lem}\label{lem:triangle}
Let $\Delta \subset \C$ be any non-degenerate triangle with sides $a,b,c \in \C$ such that $a + b + c = 0$. 
Let us denote by $e^{\sqrt{-1} \varphi / 2} \Delta$ the triangle obtained by rotating $\Delta$ by angle $\varphi / 2$ in the positive direction, with sides 
$e^{\sqrt{-1} \varphi / 2} a, e^{\sqrt{-1} \varphi / 2} b, e^{\sqrt{-1} \varphi / 2}c$. 
Then, for each side $a,b,c$ there exists exactly one value $\varphi_a, \varphi_b, \varphi_c \in [0, 2\pi )$ such that $e^{\sqrt{-1} \varphi_a / 2} a$ 
(respectively $e^{\sqrt{-1} \varphi_b / 2}b, e^{\sqrt{-1} \varphi_c / 2}c$) is purely imaginary. 
We have 
$$
  \Re ( e^{\sqrt{-1} \varphi_a / 2} b ) = \Re ( e^{\sqrt{-1} \varphi_a / 2} c ). 
$$
In addition, the function 
$$
  \Re ( e^{\sqrt{-1} \varphi_a / 2} b ) - \Re ( e^{\sqrt{-1} \varphi_a / 2} c )
$$
changes sign at $\varphi = \varphi_a$. Similar statements hold with $a,b,c$ permuted. 
\end{lem}
We call $\varphi_a, \varphi_b, \varphi_c$ the \emph{critical angle} of the sides $a,b,c$ respectively. 
By genericity, the critical angles $\varphi_a, \varphi_b, \varphi_c$ are pairwise different. 
It follows from the lemma that the critical angles decompose $S^1$ into three closed arcs 
$$
  S^1 = I_1 \cup I_2 \cup I_3 
$$
pairwise intersecting each other in a critical angle, satisfying the property: 
$$
  \max (|\Re (e^{\sqrt{-1} \varphi / 2} (\pi_0 - \pi_1 ))|, |\Re (e^{\sqrt{-1} \varphi / 2} (\pi_t - \pi_0 ))|, |\Re ( e^{\sqrt{-1} \varphi / 2}(\pi_1 - \pi_t))| )
$$
is realized
\begin{itemize}
 \item by $|\Re (e^{\sqrt{-1} \varphi / 2} (\pi_0 - \pi_1 ))|$ for $\varphi \in I_1$, 
 \item by $|\Re (e^{\sqrt{-1} \varphi / 2} (\pi_t - \pi_0 ))|$ for $\varphi \in I_2$, 
\item and by $|\Re (e^{\sqrt{-1} \varphi / 2} (\pi_1 - \pi_t ))|$ for $\varphi \in I_3$.
\end{itemize}
Namely, $I_1$ is the arc with end-points $\varphi_b, \varphi_c$ not containing $\varphi_a$, and so on. 
Let us denote by $\mbox{Int} ( I )$ the interior of an arc $I \subset S^1$, and for ease of notation let us set 
$X_j = X_j(R, \varphi)$.
It follows that  as $R \to + \infty$ 
\begin{itemize}
 \item for $\varphi \in \mbox{Int} ( I_1 )$, we have 
 \begin{align*}
  \frac{X_1}{X_2} & \to \infty,  & \frac{X_1}{X_3} & \to \infty, & 
  [X_0 : X_1 : X_2 : X_3] & \to [0:1:0:0]
 \end{align*}
\item for $\varphi \in \mbox{Int} (I_2 )$, we have 
 \begin{align*}
  \frac{X_2}{X_1} & \to \infty, & \frac{X_2}{X_3} & \to \infty, & 
  [X_0 : X_1 : X_2 : X_3] & \to [0:0:1:0],
 \end{align*}
 \item 
 for $\varphi \in \mbox{Int} ( I_3 )$, we have 
 \begin{align*}
  \frac{X_3}{X_1} & \to \infty, & \frac{X_3}{X_2} & \to \infty, &
  [X_0 : X_1 : X_2 : X_3] & \to [0:0:0:1],
 \end{align*} 
\end{itemize}
all convergence rates being exponential in $\sqrt{R}$. 
These limits show that 
\begin{itemize}
 \item for $\varphi \in \mbox{Int} ( I_1 )$, we have 
 $$
  \phi_1 \circ \psi (\E_{R, \varphi} ,\theta_{R, \varphi}) = 0, 
 $$
 \item for $\varphi \in \mbox{Int} (I_2 )$, we have 
 $$
  \phi_2 \circ \psi (\E_{R, \varphi} ,\theta_{R, \varphi}) = 0, 
 $$
 \item for $\varphi \in \mbox{Int} ( I_3 )$, we have 
 $$
  \phi_3 \circ \psi (\E_{R, \varphi} ,\theta_{R, \varphi}) = 0. 
 $$
\end{itemize}
Said differently, 
\begin{itemize}
 \item for $\varphi \in \mbox{Int} ( I_1 )$, we have 
 $$
  \phi \circ \psi (\E_{R, \varphi} ,\theta_{R, \varphi}) \in [v_2 v_3], 
 $$
 \item for $\varphi \in \mbox{Int} (I_2 )$, we have 
 $$
  \phi \circ \psi (\E_{R, \varphi} ,\theta_{R, \varphi}) \in [v_3 v_1], 
 $$
 \item for $\varphi \in \mbox{Int} ( I_3 )$, we have 
 $$
  \phi \circ \psi (\E_{R, \varphi} ,\theta_{R, \varphi}) \in [v_1 v_2].  
 $$
\end{itemize}
We infer that as 
$\varphi$ ranges over $[0,2\pi ]$ 
the corresponding elements 
$$
  \phi \circ \psi (\E_{R, \varphi} ,\theta_{R, \varphi}) \in {\Mod}_{\Betti}
$$ 
describe a path which is a generator of $\pi_1 (| \mathcal{N} (\bar{D}_{\infty}) |) \cong \Z$. 
This finishes the proof of Theorem~\ref{thm:Simpson}.  

\section{Discussion}
It follows from the proof that roughly speaking, in the large $R$ limit non-abelian Hodge theory maps the interior of each arc $I_1, I_2, I_3$ into 
a neighbourhood of one of the three intersection points of the  components $L_1, L_2, L_3$ of $\bar{D}_{\infty}$, and tiny intervals of the form 
$(\varphi_a - \varepsilon, \varphi_a + \varepsilon)$ centered at the critical angle $\varphi_a$ into one of the components (say $L_1$) of $\bar{D}_{\infty}$, 
in such a way that the image of $(\varphi_a - \varepsilon, \varphi_a + \varepsilon)$ connect $L_1 \cap L_2 = [0:0:1:0]$ with $L_1 \cap L_3 = [0:0:0:1]$. 
In different terms, large open arcs of the boundary sphere of the Hitchin base are mapped to edges of the dual boundary complex on the Betti side, 
corresponding to small neighbourhoods around the intersection points of the divisor components; 
conversely, a small open interval around a critical angle gets mapped onto the image under $\varphi$ of one of the components $L_j$. 
Such a ``large to small'' phenomenon has been conjecturally described in~\cite{KNPS} and~\cite{KS}, and is related to wall-crossing formulas.

\end{document}